\documentclass{ceb-l}
\usepackage{amssymb}

\usepackage{graphicx}

\newtheorem{theorem}{Theorem}[section]

\newtheorem{proposition}[theorem]{Proposition}

\theoremstyle{definition}
\newtheorem{definition}[theorem]{Definition}

\newtheorem{question}[theorem]{Question}
\newtheorem{conjecture}[theorem]{Conjecture}

\theoremstyle{remark}

\renewcommand{\leq}{\leqslant}
\renewcommand{\geq}{\geqslant}
\newcommand\SL{\operatorname{SL}}
\newcommand\GL{\operatorname{GL}}
\newcommand\SO{\operatorname{SO}}
\newcommand\tr{\operatorname{tr}}
\newcommand\id{\operatorname{id}}
\newcommand\poly{\operatorname{poly}}
\def\F{\mathbb{F}}
\def\R{\mathbb{R}}
\def\C{\mathbb{C}}
\def\Z{\mathbb{Z}}
\def\E{\mathbb{E}}
\def\P{\mathbb{P}}
\def\Q{\mathbb{Q}}

\numberwithin{equation}{section}

\begin{document}

\title[Approximate groups and their applications]{Approximate groups and their applications: work of Bourgain, Gamburd, Helfgott and Sarnak}


\author{Ben Green}
\address{Centre for Mathematical Sciences, Wilberforce Road, Cambridge CB3 0WA}
\curraddr{Radcliffe Institute for Advanced Study, 8 Garden Street, Cambridge MA 02138}
\email{b.j.green@dpmms.cam.ac.uk}
\thanks{This article was written while the author was a fellow at the Radcliffe Institute at Harvard. It is a pleasure to thank the institute for its support and excellent working conditions.}

\subjclass[2000]{Primary }

\begin{abstract}
This is a survey of several exciting recent results in which techniques originating in the area known as additive combinatorics have been applied to give results in other areas, such as group theory, number theory and theoretical computer 
science.  We begin with a discussion of the notion of an approximate group and also that of an approximate field, describing key 
results of Fre\u{\i}man-Ruzsa, Bourgain-Katz-Tao, Helfgott and others in which the structure of such objects is elucidated. 
We then move on to the applications. In particular we will look at the work of Bourgain and Gamburd on expansion 
properties of Cayley graphs on $\mbox{SL}_2(\F_p)$ and at its application in the work of Bourgain, Gamburd and Sarnak on nonlinear sieving problems. 
\end{abstract}

\maketitle

\section{Introduction}

The subject of \emph{additive combinatorics} has grown enormously over the last ten years and now comprises a large collection of tools with many applications in number theory and elsewhere, for example in group theory and theoretical computer science. It has often been thought a little difficult to specify to an outsider exactly what the subject \emph{is}\footnote{See, for example, my own attempt in the opening remarks of \cite{green-review}.}. However the following point of view seems to be gradually crystallising: additive combinatorics is the study of \emph{approximate mathematical structures} such as approximate groups, rings, fields, polynomials and homomorphisms. It is interested in what the right definitions of these approximate structures are, what can be said about them, and what applications this has to other parts of mathematics. 

This article has three main aims. Firstly, we wish to introduce the above point of view to a general audience, focussing in particular on the basic theory of approximate groups and approximate fields. Secondly, we wish to sketch some beautiful applications of these ideas. One of them has to do with the beautiful picture on the cover (for which we thank Cliff Reiter) of an Apollonian circle packing. It is classical that the radii of these circles are all reciprocals of integers. We will describe work of Bourgain, Gamburd and Sarnak giving upper bounds for the number of circles at ``depth $n$'' which have radius the reciprocal of a prime. Thirdly, we wish to hint at the extraordinary variety of different areas of mathematics which have started to interact with additive combinatorics: geometric group theory, analytic number theory, model theory and point-set topology are just the ones we shall mention here. 

What we offer here is merely a taste of this viewpoint of additive combinatorics as the theory of approximate structure and of its applications. We do not touch on the theory of approximate polynomials (a.k.a. the theory of Gowers norms) or say much at all about approximate homomorphisms, or anything about the many applications of these two notions. These topics will be covered in detail in forthcoming lecture notes of the author \cite{green-barbados}.

\section{Approximate groups}\label{approx-groups-sec}

Before we can define an approximate group, we need to recall what an exact one is. We shall be concerned with finite groups, and we shall be working inside some ambient group $G$, so that it makes sense to talk about multiplication of elements and taking inverses. If $A \subseteq G$ is a finite set then we shall write $A \cdot A := \{a_1 a_2 : a_1,a_2 \in A\}$ and $A \cdot A^{-1} := \{a_1 a_2^{-1} : a_1,a_2 \in A\}$. Later on we shall see more general notations such as $A \cdot A \cdot A$ and $A \cdot B$ whose meaning, we hope, will be evident. The following proposition, whose proof is an exercise in undergraduate group theory, gives various criteria for $A$ to be a subgroup or something very closely related.

\begin{proposition}\label{subgroup-char} Let $A$ be a finite subset of some ambient group $G$. Then we have the following statements\footnote{Statements (5) and (6) look ``probabilistic'' but this is just a notation. By $\P(a_1a_2 \in A| a_1, a_2 \in A)$ we mean simply the proportion of all pairs $a_1,a_2 \in A$ for which $a_1 a_2$ also lies in $A$.}:
\begin{enumerate}
\item $|A \cdot A^{-1}| \geq |A|$, with equality if and only if $A = Hx$ for some subgroup $H \leq G$ and some element $x \in G$;
\item $|A \cdot A| \geq |A|$, with equality if and only if $A = Hx$ for some subgroup $H \leq G$ and some element $x$ in the normaliser $N_G(H)$;
\item The number of quadruples $(a_1,a_2,a_3,a_4) \in A^4$ with $a_1 a_2^{-1} = a_3 a_4^{-1}$ is at most $|A|^3$, with equality if and only if $A = Hx$ for some subgroup $H \leq G$ and some element $x \in G$;
\item The number of quadruples $(a_1,a_2,a_3,a_4) \in A^4$ with $a_1 a_2 = a_3 a_4$ is at most $|A|^3$, with equality if and only if $A = Hx$ for some subgroup $H \leq G$ and some element $x \in N_G(H)$;
\item $\P(a_1 a_2 \in A | a_1, a_2 \in A) \leq 1$, with equality if and only if $A = H$ for some subgroup $H \leq G$;
\item $\P(a_1 a_2^{-1} \in A | a_1, a_2 \in A) \leq 1$, with equality if and only if $A = H$ for some subgroup $H \leq G$.
\end{enumerate}
\end{proposition}

This would be a rather odd proposition to see formulated in an algebra text. However each of the statements (1) -- (6) has been constructed as an inequality in such a way that one may ask when equality approximately holds. Before we can talk about such approximate variants, however, we need to know how approximate they will be. For this purpose we introduce a parameter $K \geq 1$; larger values of $K$ will indicate more approximate, and thus less structured, objects\footnote{In practice the theory when $K \approx 1$ is very different from the theory when, for example, $K \sim 100$. In the former setting, these approximate notions of subgroup constitute very small perturbations of the exact characterisations of Proposition \ref{subgroup-char}, and it turns out (though is not always trivial to prove) that the approximate objects so defined are small perturbations of the exact objects characterised by Proposition \ref{subgroup-char}. In conversation Tao and I tend to refer to this regime as ``the 99\% world'', an expression I would not be averse to popularising. In this paper $K$ will be much larger, causing the theory to become much richer. Tao and I call this the ``1\% world'' although the parameter $K$ could be anything between $2$ (say) and some small power of $|A|$.}. 

Consider, then, the following list of properties that a finite set $A \subseteq G$ might enjoy.

\begin{enumerate}
\item $|A \cdot A^{-1}| \leq K |A|$;
\item $|A \cdot A| \leq K |A|$;
\item The number of quadruples $(a_1,a_2,a_3,a_4) \in A^4$ with $a_1 a_2^{-1} = a_3 a_4^{-1}$ is at least $|A|^3/K$;
\item The number of quadruples $(a_1,a_2,a_3,a_4) \in A^4$ with $a_1 a_2 = a_3 a_4$ is at least $|A|^3/K$;
\item $\P(a_1 a_2 \in A | a_1, a_2 \in A) \geq 1/K$;
\item $\P(a_1 a_2^{-1} \in A | a_1, a_2 \in A) \geq 1/K$.
\end{enumerate}

Now these are by no means as closely equivalent as the properties (1) -- (6) in Proposition \ref{subgroup-char}. Let us give an example in which the ambient group is $\Z$, and where we use additive rather than multiplicative notation. Take $A = \{1,\dots,n\} \cup \{2^{n+1}, 2^{n+2},\dots,2^{2n}\}$. Then it is easy to check that (3) and (4) are both satisfied with any $K > 12$, as $n \rightarrow \infty$, this being because there are $\frac{2}{3}n^3 (1 + o(1))$ solutions to $a_1 + a_2 = a_3 + a_4$ with $a_1,a_2,a_3,a_4 \in \{1,\dots,n\}$. On the other hand the sumset $A + A$ contains the numbers $2^{n + i} + j$ for each pair $i,j$ with $0 < i,j \leq n$. Since these numbers are all distinct, we have $|A + A| \geq n^2 = |A|^2/2$, which means that if $n$ is sufficiently large depending on $K$ then (2) is not satisfied at all.

Rather remarkably, however, there is a sense in which the concepts (1) -- (6) are all \emph{roughly} the same. To say what we mean by that, we introduce the following notion of rough equivalence\footnote{The fact that we have written $Bx$ rather than $xB$ is a little arbitrary. The notion of rough equivalence will, in this survey, be applied to classes of sets (such as (1) -- (6) here) which are invariant under conjugation, in which case whether we multiply on the left or on the right in the definition makes little difference.}.

\begin{definition}[Rough Equivalence]
Suppose that $A$ and $B$ are two finite sets in some ambient group and that $K \geq 1$ is a parameter. Then we write $A \sim_K B$ to mean that there is some $x$ in the ambient group such that $|A \cap Bx| \geq \max(|A|,|B|)/K$. We say that $A$ and $B$ are roughly equivalent \textup{(}with parameter $K$\textup{)}.
\end{definition}

The remarkable fact alluded to above is the following. For every choice of $j,j' \in \{1,\dots,6\}$, suppose that some set $A$ satisfies condition ($j$) in the list above with parameter $K$. Then there is a set $B$ satisfying condition ($j'$) with parameter $K' = \poly(K)$ (some polynomial in $K$) such that $A \sim_{K'} B$.  Of particular note is the fact that the weak ``statistical'' properties (3) -- (6) imply the apparently more structured properties (1) and (2). The proof of this is not at all trivial and the main content of it is the so-called Balog-Szemer\'edi-Gowers theorem \cite{gowers1}, generalised to the nonabelian setting in the fundamental paper of Tao \cite{tao-noncommutative}, as well as a collection of ``sumset estimates'' which, in the abelian case, I refer to collectively as \emph{Ruzsa calculus} \cite{green-partiii}. These estimates of Ruzsa have such a classical role in the theory that we record two of them, in the abelian setting, here: we will mention these two again later on. 
\begin{theorem}[Ruzsa]\label{ruzsa-sumsets}
Suppose that $A_1,A_2$ and $A_3$ are finite sets in some ambient abelian group. Then $|A_1||A_2-A_3| \leq |A_1 - A_2||A_1 - A_3|$ and $|A_1 + A_1| \leq |A_1-A_1|^3/|A_1|^2$.
\end{theorem}
The original paper \cite{ruzsa-sumset-1}, the book \cite{tao-vu-book} or the notes \cite{green-partiii} may be consulted for more details. The first estimate is true in general groups but adapting the second requires care: see \cite{tao-noncommutative}. 

It might be remarked that for many pairs ($j$) and ($j'$) the correspondence between the relevant properties is a little tighter than mere rough equivalence, and often this can be useful. We shall not dwell on this point here. In the paper of Tao just mentioned one finds what has become the ``standard'' notion of an approximate group.

\begin{definition}[Approximate group]\label{approx-gp-def}
Suppose that $A$ is a finite subset of some ambient group and that $K \geq 1$ is a parameter. Then we say that $A$ is a $K$-approximate group if it is symmetric (that is, if $a \in A$ then $a^{-1} \in A$, and the identity lies in $A$) and if there is a set $X$ in the ambient group with $|X| \leq K$ and such that $A\cdot A \subseteq X \cdot A$. 
\end{definition}

This notion, it turns out, is roughly equivalent to (1) - (6) above. It has certain advantages over (1) - (6), for example as regards its behaviour under homomorphisms. It is also clear that an approximate group in this sense enjoys good control of iterated sumsets. Thus, for example, $A \cdot A \cdot A \subseteq X \cdot X \cdot A$, which means that $|A^3| = |A \cdot A \cdot A| \leq K^2 |A|$, and similarly $|A^n| \leq K^{n-1}|A|$ where $A^n$ denotes the set of all products $a_1\dots a_n$ with $a_1,\dots,a_n \in A$. From now on, when we speak of an approximate group, we will be referring primarily to Definition \ref{approx-gp-def}.

With this discussion in mind, we can introduce what might be termed the rough classification problem of approximate group theory.

\begin{question}\label{rough-classification-problem}
Consider the collection $\mathcal{C}$ of all $K$-approximate groups $A$ in some ambient group $G$. Is there some ``highly structured'' subcollection $\mathcal{C}'$ such that every $A \in \mathcal{C}$ is roughly equivalent to some set $B \in \mathcal{C}'$ with parameter $K'$, where $K'$ depends only on $K$?
\end{question}

This question has been addressed in a great many different contexts, starting with the Fre\u{\i}man-Ruzsa theorem \cite{freiman-book,ruz-frei}, which gives an answer for subsets of $\Z$. Here, it is possible to take $\mathcal{C}'$ to consist of the so-called \emph{generalised arithmetic progressions}, that is to say sets $B$ of the form
\[ B := \{ l_1 x_1 + \dots + l_d x_d : l_i \in \Z, |l_i| \leq L_i\},\]
where $x_1,\dots,x_d \in \R$, the quantities $L_1,\dots,L_d$ are ``lengths'' and $d \leq K$. Note in particular that, even in the highly abelian setting of the integers $\Z$, approximate groups are a more general kind of object than genuine subgroups. That is, the theory of approximate groups, even up to rough equivalence, is a little richer than the theory of finite subgroups of $\Z$ (which is in fact a rather trivial theory). The remarkable feature of the Fre\u{\i}man-Ruzsa theorem is that the theory is not \emph{much} richer, in the sense that generalised progressions remain highly ``algebraic'' objects. Here is a list of other contexts in which the question has been at least partially answered: 
\begin{itemize}
\item abelian groups \cite{green-ruzsa};
\item nilpotent and solvable groups \cite{bg1,bg2,fisher-katz-tao,sanders,tao-solvable};
\item free groups \cite{razborov};
\item linear groups $\mbox{SL}_2(\R)$ \cite{elekes-kiraly},\\
 $\mbox{SL}_2(\C)$ \cite{chang-sl3,helfgott-1,helfgott-2}, \\
 $\SL_3(\Z)$ \cite{chang-sl3},\\
 $\SL_3(\C)$ (sketched in \cite{helfgott-2}),\\
 ``bounded'' subsets of $\SL_n(\C)$ including $\mbox{U}_n(\C)$ \cite{bg3},\\
  $\mbox{SL}_2(\F_p)$ \cite{helfgott-1}, \\
  $\mbox{SL}_3(\F_p)$ \cite{helfgott-2}\\
  and $\SL_2(\Z/q\Z)$  for various other $q$ (cf. \cite{bourgain-gamburd-pn}). 
\end{itemize}

It is generally felt that approximate groups in quite general contexts can be controlled by objects built up from genuine subgroups and nilpotent objects; this has been found in all of the examples just mentioned and is suggested by the famous theorem of Gromov on groups with polynomial growth \cite{gromov} and the recent quantitative formulation of it due to Shalom and Tao \cite{shalom-tao}. Quite precise suggestions along these lines have been made by Helfgott, Lindenstrauss and others: more information on this can be found on Tao's blog \cite{tao-blog}.

Before leaving this subject, we remark that even (perhaps \emph{especially}) in the abelian case the issue of the dependence of $K'$ on $K$ is far from being resolved. No examples are known to rule out the possibility that, with the right definition of the highly-structured class $\mathcal{C}'$, $K'$ can be taken to be polynomial in $K$. In particular this is conjectured when the ambient group $G$ is $\F_2^{\Z}$, the countable infinite vector space over the field of two elements, and $\mathcal{C}'$ consists of (finite) subgroups of $G$. This assertion\footnote{There are variants of this conjecture over other groups, such as $\Z$; see \cite{gowers-gafa,green-tao-equivalence}.} is known as the \emph{polynomial Fre\u{\i}man-Ruzsa conjecture} \cite{ruzsa-finite-field}, see also \cite{green-finite-field,pfr-blog}. It is equivalent to the following question which, for many years, I have tried to advertise to those for whom the word cohomology holds no fear.

\begin{question}
Suppose that $\phi : \F_2^n \rightarrow \F_2^{\Z}$ is a map such that $\phi(x+y) - \phi(x) - \phi(y)$ takes on at most $K$ different values as $x,y$ range over $\F_2^n$. Is it true that $\phi = \tilde \phi + \eta$, where $\tilde \phi$ is linear and $\eta$ takes on at most $K' = \poly(K)$ different values?
\end{question}

It is a very easy exercise to obtain such a statement with $K' = 2^K$ but, so far as I know, no serious improvement of this bound has ever been obtained\footnote{I would be very interested to see even a bound of the form $2^{o(K)}$.}.

\section{Approximate rings and fields}

Fortified by the experiences of the last section, one might attempt to come up with a sensible notion of an approximate \emph{ring}. A natural one, based perhaps on (2) in the previous section, is as follows: if $A$ is a finite subset of some ambient ring $R$, we say that it is a $K$-approximate ring if $|A + A| \leq K|A|$ and $|A \cdot A| \leq K|A|$. Here, of course, $A + A := \{a_1 + a_2 : a_1,a_2 \in A\}$ and $A \cdot A = \{a_1 a_2 : a_1,a_2 \in A\}$ as before. If $R = \F$ is actually a field (or an integral domain, which embeds into its field of fractions) then we refer to $A$ as an approximate field, noting that approximate closure under division is essentially automatic in view of the rough equivalence of the notions (1) and (2) of approximate group.

The study of approximate rings and fields was initiated in a paper of Erd\H{o}s and Szemer\'edi \cite{erdos-szemeredi} who proved (though not in this language!) that a $K$-approximate subfield of $\Z$ must have size $\poly(K)$. They in fact conjectured that the right bound is $C_{\epsilon} K^{1 + \epsilon}$ for any $\epsilon > 0$, but this is so far unresolved; the best exponent so far obtained is $3 + \epsilon$, a result of Solymosi \cite{solymosi-43}. Note that this is equivalent to, and more usually stated as, the lower bound 
\[ \max(|A + A|, |A\cdot A|) \geq c_{\epsilon} |A|^{4/3 + \epsilon}\] for all finite sets $A \subseteq \Z$. In a different paper \cite{solymosi-c}, Solymosi generalised the Erd\H{o}s-Szemer\'edi result to show that every $K$-approximate subfield of $\C$ has size at most $2^{12}K^4$. 

The general theory of approximate fields can be said to have started with the papers of Bourgain-Katz-Tao \cite{bkt} and Bourgain-Glibichuk-Konyagin \cite{bourgain-glibichuk-konyagin,bourgain-konyagin}, where\footnote{The original paper \cite{bkt} of Bourgain, Katz and Tao did not quite classify the very small (smaller than $p^{\delta}$) approximate subrings of $\F_p$; this restriction was removed in \cite{bourgain-glibichuk-konyagin,bourgain-konyagin}. Very often the approximate fields under consideration in a given setting will have size at least $p^{\delta}$, and for this reason one often refers to the Bourgain-Katz-Tao theorem.}  the following result is established.

\begin{theorem}\label{bkt-theorem}
Let $p$ be a prime and let $K \geq 2$. Then every $K$-approximate subfield of $\F_p$ has size at most $K^C$ or at least $K^{-C}p$, for some absolute constant $C$.
\end{theorem}

The arguments on page 384 of \cite{bourgain-glibichuk-konyagin}, though they are phrased in a more limited context, essentially prove that every approximate subfield (in an arbitrary ambient field) must be roughly equivalent to a genuine finite subfield. This unifies the results of Erd\H{o}s-Szemer\'edi and Solymosi with Theorem \ref{bkt-theorem}. In fact something similar is true for approximate rings, at least provided the ambient ring $R$ does not have ``too many'' zero divisors. These issues are comprehensively explored in an interesting paper \cite{tao-general-rings} of Tao, which also has a very comprehensive collection of references.

Suppose that $A$ is a $K$-approximate field in some ambient field $\F$, that is to say both $|A \cdot A|$ and $|A + A|$ are bounded by $K|A|$. We are going to sketch a proof that $\F$ must contain a genuine subfield $B$ which is ``close'' to $A$. The first step is to prove the \emph{Katz-Tao lemma}, which asserts that $A$ (or, more precisely, a large subset $A' \subseteq A$) behaves in a manner which more strongly resembles that of a field: that is to say, $A$ is almost closed under both addition/subtraction and multiplication/division \emph{simultaneously}. To give a (relevant) example, the set
\[ \overline{A} := \{ a_5 + a_6\frac{a_1 - a_3}{a_4 - a_2} : a_1,\dots,a_6 \in A\}\] has size $\overline{K}|A|$, where $\overline{K} = \poly(K)$. 

A slick proof of the Katz-Tao lemma is given in \cite[Section 2.5]{tao-general-rings} and we shall say little more about it here other than to remark that it involves a combination of Ruzsa's sumset calculus and clever elementary arguments. Personally, I regard it as part of the ``basic'' theory of approximate fields as opposed to the ``structural theory'', to be regarded on the same level as the arguments used to show that definitions (1) -- (6) of an approximate group are roughly equivalent (namely, Ruzsa's sumset calculus and the Balog-Szemer\'edi-Gowers theorem). In other words one might argue that the smallness of $\overline{A}$, or of similar objects, might be taken as an alternative \emph{definition} of approximate field.

Suppose that $A$ is known to have this property, that is to say $|\overline{A}| \leq \overline{K}|A|$. Then it is possible to establish an intriguing dichotomy: if $\xi \in \F^{\times}$ then either
\begin{equation}\label{poss-1} |A + A\xi| = |A|^2\end{equation}
or 
\begin{equation}\label{poss-2} |A + A\xi| \leq \overline{K}|A|.\end{equation}
Here, $A + A\xi$ refers to the set of all $a_1 + a_2\xi$ with $a_1,a_2 \in A$. To see why this is so, note that $|A + A\xi| \leq |A|^2$ and that equality occurs if and only if the elements $a_1 + a_2\xi$ are all distinct. If equality does not occur then we may find a nontrivial solution to $a_1 + a_2\xi = a_3 + a_4\xi$, which means that $\xi = \frac{a_1 - a_3}{a_4 - a_2}$. But then every element of $A + A\xi$ has the form 
\[ a_5 + a_6\xi = a_5 +  a_6\frac{a_1 - a_3}{a_4 - a_2} ,\] and thus lies in $\overline{A}$. 

On the other hand, it is not hard to see using Ruzsa calculus\footnote{In addition to the bounds of Theorem \ref{ruzsa-sumsets} one requires an inequality controlling $|A_1 + A_2 + A_3|$ in terms of the $|A_i + A_j|$.} that if $\xi_1,\xi_2$ satisfy \eqref{poss-2} then for $\xi = \xi_1+ \xi_2, \xi_1 - \xi_2,\xi_1\xi_2,\xi_1 \xi^{-1}_2$ we have \[ |A + \xi \cdot A| \leq \overline{K}^C |A| \leq K^{C'}|A|\] for absolute constants $C,C'$. If $K$ is a sufficiently small power of $|A|$ then this means that \eqref{poss-1} cannot hold, forcing us to conclude that \eqref{poss-2} holds for $\xi$. In this way we identify the set\footnote{Note that this set may be identified with $\frac{A - A}{(A - A)^{\times}}$.} of all $\xi$ satisfying \eqref{poss-2} as a genuine subfield of $\F$. Straightforward additional arguments allow one to show that this subfield and $\F$ are roughly equivalent.\vspace{11pt}

The original argument of \cite{bkt} is different and specific to $\F_p$ but rather fun and, given the preceding discussion, it is not hard to say a few meaningful words about it.  Suppose for the sake of illustration that $A \subseteq \F_p$ is a $K$-approximate subfield of size $\sim p^{1/10}$; our task is to derive a contradiction if (say) $K = p^{o(1)}$. Suppose that the Katz-Tao lemma has already been applied, so that $\overline{A}$, as defined above, is known to be small. The sets $\overline{\overline{A}}, \overline{\overline{\overline{A}}},\dots$ arising from (boundedly many) successive applications of this operation may also be shown to be small. Now simple averaging arguments (using nothing more than the fact that $|A| = p^{1/10}$) show that $\F_p$ has dimension at most 100 (say) as a ``vector space'' over $A$; that is, there exist $x_1,\dots,x_{100} \in \F_p$ such that \begin{equation}\label{span}\F_p = Ax_1 + \dots + Ax_{100}.\end{equation} Now $x_1,\dots,x_{100}$ cannot be a ``basis'' for $\F_p$ over $A$ since otherwise we would have $p = |A|^{100}$, contrary to the assumption that $p$ is prime. Thus there must exist some $x \in \F_p$ which is representable in two different ways as
\[ x = a_1 x_1 + \dots + a_{100} x_{100} = a'_1 x_1 + \dots + a'_{100} x_{100}\] with $a_1,\dots,a_{100}, a'_1,\dots,a'_{100} \in A$.
Suppose, without loss of generality, that $a_{100} \neq a'_{100}$. Then 
\[ x_{100} = \frac{(a_1 - a'_1)x_1 + \dots + (a_{99} - a'_{99})x_{99}}{a'_{100} - a_{100}}.\]
By substituting this expression for $x_{100}$ into \eqref{span}, we see that 
\[ \F_p = \overline{A} x_1 + \dots + \overline{A} x_{99}.\] 
Repeating the argument gives (without loss of generality) \[ \F_p = \overline{\overline{A}} x_1 + \dots + \overline{\overline{A}} x_{98},\] and we may continue in this fashion to get, eventually,
\[ \F_p = \overline{\overline{\dot{\dot{\overline{\overline{A}}}}}} x_1.\]
This is contrary to the fact that none of the sets $\overline{A},\overline{\overline{A}},\dots$ has size much larger than that of $A$ itself, namely about $p^{1/10}$, and a contradiction ensues.\vspace{11pt}

Remarkably, the main ``dimension reduction'' idea here comes from a paper in point-set topology, namely Edgar and Miller's solution of the Erd\H{o}s-Volkmann ring problem \cite{edgar-miller} (that is, the statement that all Borel subrings of $\R$ have dimension $0$ or $1$). See in particular Lemma 1.3 of that paper.

\section{Helfgott's results}

In this section we discuss the results of Helfgott \cite{helfgott-1,helfgott-2} concerning approximate subgroups of \[ \mbox{SL}_2(\F_p) := \{\begin{pmatrix} a & b \\ c & d \end{pmatrix} : a,b,c,d \in \F_p : ad - bc = 1\}.\] Helfgott proves the following.

\begin{theorem}[Helfgott]\label{helfgott-theorem}
Suppose that $A \subseteq \SL_2(\F_p)$ is a $K$-approximate group. Then $A$ is roughly $K^C$-equivalent to an upper-triangular $K^C$-approximate subgroup of $\SL_2(\F_p)$ \textup{(}that is, an approximate subgroup conjugate to a set of upper-triangular matrices\textup{)}. 
\end{theorem}

Rather than discuss Helfgott's result itself, we discuss the analogous question for $\SL_2(\C)$. Here the answer is rather simpler and is given in \cite{chang-sl3}, based on Helfgott's work.

\begin{theorem}\label{sl2-c-theorem}
Suppose that $A \subseteq \SL_2(\C)$ is a $K$-approximate group. Then $A$ is roughly $K^C$-equivalent to an \emph{abelian} $K^C$-approximate subgroup  of $\SL_2(\C)$. 
\end{theorem}

If desired the abelian approximate group could itself be controlled by a generalised progression using the Fre\u{\i}man-Ruzsa theorem.

We will only sketch a proof of the weaker result that $A$ is $K^C$-equivalent to an upper-triangular $K^C$-approximate subgroup, that is to say the direct analogue of Helfgott's result. In $\SL_2(\C)$, additional arguments may then be applied to prove Theorem \ref{sl2-c-theorem}; there are no such arguments in $\SL_2(\F_p)$, since the upper-triangular ``Borel subgroup'' 
\[\{\begin{pmatrix} \lambda & \mu \\ 0 & \lambda^{-1}\end{pmatrix} : \lambda \in \F_p^*, \mu \in \F_p\}\] is not close to abelian.

The proof of this weak form of Theorem \ref{sl2-c-theorem} is simpler than that of Theorem \ref{helfgott-theorem} in two major ways. Firstly since $\C$ is algebraically closed we may talk about eigenvalues, eigenvectors and diagonalization without the need to pass to an extension field, whereas over $\F_p$ we would have to involve the quadratic extension $\F_{p^2}$. Secondly, the structure of $K$-approximate subfields of $\C$ is easy to describe: by the theorem of Solymosi \cite{solymosi-c} they are all sets of size at most $2^{12}K^4$. Theorem \ref{bkt-theorem}, by contrast, has to allow for those approximate fields which are almost all of $\F_p$. Worse still, to handle $\SL_2(\F_p)$ Helfgott must in fact classify approximate subfields of $\F_{p^2}$, and this involves the additional possibility of sets which are close to the subfield $\F_p$. 

For the sake of exposition, we will assume in the first instance that $A$ is a genuine finite subgroup of $\SL_2(\C)$; our task is to show that $A$ contains a large upper-triangular subgroup. When we have sketched how Helfgott's argument looks in this case we will remark on the additional technicalities required to make the argument ``robust'' enough to apply to $K$-approximate groups.

The key idea in Helfgott's argument, referred to by subsequent authors as \emph{trace amplification}, involves examining the set of traces
\[ \tr A := \{ \tr a : a \in A\}.\]
We will sketch a proof that a large subset of this set of traces is a $2^{24}$-approximate subfield of $\C$ of size greater than $2^{108}$. This contradicts Solymosi's theorem \cite{solymosi-c} and so we must be in one of those degenerate situations. Careful analysis of each of them leads to the conclusion that $A$ is roughly upper-triangular.

The first degenerate situation to analyse is that in which $\tr A$ is small, an appropriate notion of \emph{small} being $|\tr A| \leq 2^{111}$. Now a linear algebra computation (Lemma 4.2 of \cite{bourgain-gamburd-sarnak}) confirms that if $g,h \in A$ are elements without a common eigenvector in $\C^2$ then the map \[ \SL_2(\C) \rightarrow \C^3 : x \mapsto (\tr x, \tr(gx), \tr(hx))\] is at most two-to-one. This, or rather the fact that something like this holds, is not at all surprising: indeed knowledge of $\tr(x), \tr(gx), \tr(hx)$ together with the fact that $\det(x) = 1$ provides four pieces of information which, generically, ought to more-or-less determine the four entries of the matrix $x$. If $A$ contains two such elements $g,h$ then it follows that we have 
\[|A| \leq 2|\tr A|^3 \leq 2^{334},\] and so $|A|$ is also small\footnote{Additive combinatorics has a bad reputation for referring to quantities like $2^{334}$ as ``small''. ``Bounded by an absolute constant'' might be more appropriate.}. If, by contrast, $A$ does not contain two such elements, and if $|A| > 3$, then it is easy to see that there is some $v \in \C^2$ which is an eigenvector for all of $A$ simultaneously. With respect to a basis containing $v$, every matrix in $A$ is upper-triangular. 

Suppose, then, that $|\tr A| > 2^{111}$. In particular (!) there is some element $g \in A$ which is non-parabolic, or in other words $\tr g \neq \pm 2$; such elements have distinct eigenvalues and so are diagonalisable. 


Write $A' \subseteq A$ for the set of non-parabolic elements; then $|\tr A'| \geq |\tr A| - 2 \geq \frac{1}{2}|\tr A|$. Now in $\SL_2(\C)$ the trace of a non-parabolic element $g$ completely determines the conjugacy class of $g$. It follows that there is some non-parabolic $g \in A$ such that the conjugacy class of $A$ containing $g$ has size at most $2|A|/|\tr A|$. By the orbit-stabiliser theorem, the centraliser\footnote{$T$ is for \emph{torus}, the word used for such a subgroup in Lie theory.} \[ T = C_A(g) = \{a \in A: ag = ga\}\] has size at least $\frac{1}{2}|\tr A|$. But by changing basis so that $g$ is in diagonal form (with distinct diagonal entries) it is not hard to check that $T$ consists entirely of diagonal matrices. No single trace can arise from more than two of these elements, and so $|\tr T| \geq \frac{1}{4}|\tr A| > 2^{109}$. We shall show that the set
\[ R := \{ \tr a^2 : a \in T\}\]
is a $2^{24}$-approximate subfield of $\C$. Noting that
\begin{equation}\label{trt} |R| \geq \textstyle\frac{1}{2}\displaystyle|\tr T| > 2^{108},\end{equation}this is contrary to Solymosi's theorem. In order to do this we play around a little with traces. Such playing around is most productive if, in the basis just selected, there is an element $a = \begin{pmatrix} a_{11} & a_{12} \\ a_{21} & a_{22} \end{pmatrix} \in A$ with $a_{11}a_{12}a_{21}a_{22} \neq 0$. The absence of such an element is another degenerate situation to analyse, and once again one can check\footnote{This is, admittedly, a somewhat tedious check.} that $A$ must be either upper-triangular or else equal to one of the dihedral groups, each of which has an index two abelian subgroup. 

Now let us note that 
\begin{equation}\label{ar-1} R\cdot R  \subseteq R + R,\end{equation}
this being a consequence of the fact that 
\[ (t^2_1 + t_1^{-2})(t^2_2 + t_2^{-2}) = (t^2_1 t^2_2 + t_1^{-2} t_2^{-2})  + (t^2_1 t_2^{-2} + t_1^{-2} t^2_2).\] Let us also note that
\[ \tr \big(  \begin{pmatrix} t_1t_2 & 0 \\ 0 & t_1^{-1}t_2^{-1}\end{pmatrix}  a  \begin{pmatrix} t_1t^{-1}_2 & 0 \\ 0 & t_1^{-1}t_2\end{pmatrix} a^{-1} \big) = \mu(t_1^2 + t_1^{-2}) + \lambda (t_2^2 + t_2^{-2}),\] 
where $\mu := a_{11}a_{22} \neq 0$ and $\lambda := -a_{12}a_{21} \neq 0$, which means that 
\[ \lambda R + \mu R \subseteq \tr A.\]
In particular 
\[ |R + \frac{\mu}{\lambda} R| = |\lambda R+ \mu R| \leq |\tr A| \leq 16 |R|,\] which, by Ruzsa's inequalities (Theorem \ref{ruzsa-sumsets}, applied with $A_1 =\frac{\mu}{\lambda} R$ and $A_2 = A_3 = -R$) implies that $|R + R| \leq 2^{24}|R|$. This, together with \eqref{ar-1}, implies that $R$ is a $2^{24}$-approximate subring of $\C$. By Solymosi's theorem this implies that $|R| \leq 2^{108}$, contrary to \eqref{trt}.\hfill $\Box$

In the above sketch we assumed, of course, that $A$ was actually a finite subgroup. However the argument was of a type that can be made to work for $K$-approximate groups also. To explain what we mean by this let us remark, rather vaguely, on how one or two of the steps adapt and then offer some general remarks. 

\emph{Orbit-Stabiliser theorem.} If $A$ is a group and if $x \in A$ then we used the fact that the size of the conjugacy class $\Sigma(x)$ containing $x$ and that of the centraliser $C_A(x)$ are related by $|\Sigma(x)||C_A(x)| = |A|$. In fact we only used the inequality $|C_A(x)| \geq |A|/|\Sigma(x)|$, giving us an element with large centraliser, and here is a simple way of seeing why this holds: all of the conjugates $axa^{-1}$, $a \in A$, lie in $\Sigma(x)$, and so by the pigeonhole principle there must be distinct elements $a_1,\dots,a_k \in A$, $k \geq |A|/|\Sigma(x)|$, with $a_1 x a_1^{-1} = \dots = a_k x a_k^{-1}$. But then the elements $a_i^{-1} a_1$, $i = 1,\dots,k$, centralise $x$. Now if $A$ is only a $K$-approximate group then this argument does not quite work, as there is no well-defined notion of conjugacy class. However a similar pigeonhole argument nonetheless gives us an element with large centraliser, since the conjugates $a x a^{-1}$ are all constrained to lie in $A^3$, a set of size at most $K^2|A|$. 

\emph{Escape from subvarieties.} A more interesting point concerns the location of an element of $A$ which, in a given basis, has no zero entries. Whilst this might not be \emph{a priori} possible if $A$ is only an approximate group, it \emph{is} possible to find such an element in $A^n$ for some bounded $n$ (independent of the approximation parameter $K$), and this is good enough for Helfgott's purposes. This is a special case of a nice lemma of Eskin, Mozes and Oh \cite{eskin-mozes-oh} called ``escape from subvarieties''. The point is that the group $\langle A \rangle$ generated by $A$, if it is not almost upper-triangular, contains an element with no zero entries -- indeed this fact was used in the above sketch. In other words, $\langle A \rangle$ is not contained in the subvariety of $\SL_2(\C)$ defined by
\[ V := \{\begin{pmatrix} x_{11} & x_{12} \\ x_{21} & x_{22} \end{pmatrix} : x_{11}x_{12} x_{21} x_{22} = 0\}.\]
The Eskin-Mozes-Oh result states that in such a situation we can find ``evidence'' for the non-containment of $\langle A \rangle$ inside $V$ by taking just a bounded number, depending only on $V$, of products of $A$.

It seems, then, that certain types of argument -- in some sense those involving ``bounded length'' computations in the ambient group -- adapt very well from the traditional group theory setting to approximate groups. At the moment we do not have anything approaching a precise formulation of this principle and indeed at present the passage from the ``exact'' to the approximate is as much an art as a science. Nonetheless, there seems to be merit in looking for ``bounded length'' proofs in traditional group theory which might be adapted to the approximate setting. Perhaps this is as good a place as any to mention the remarkable recent paper of Hrushovski \cite{hrushovski} in which tools from model theory have been applied to the study of approximate groups. The ramifications of that paper are not yet completely clear, but it looks as though Theorem 1.3 of that paper together with some structure theory of algebraic groups ought to lead, without too much difficulty, to a proof of the following statement.

\begin{conjecture}
Suppose that $A \subseteq \SL_n(\C)$ is a $K$-approximate group. Then there is a $K'$-approximate group $B$ which is \emph{nilpotent} and $K'$-controls $A$, where $K'$ depends only on $K$.
\end{conjecture}

It seems reasonable to conjecture that $K'$ can be taken to depend polynomially on $K$, although in their present form Hrushovski's techniques will not give this.

\section{Cayley graphs on $\mbox{SL}_2(\F_p)$}

We move on now to applications of the theory of approximate groups. In this section we discuss the paper \cite{bourgain-gamburd-sl2} of Bourgain and Gamburd. This paper concerns \emph{expander graphs}. For the purposes of this discussion these are $2k$-regular graphs $\Gamma$ on $n$ vertices for which there is a constant $c > 0$ such that for any set $X$ of at most $n/2$ vertices of $\Gamma$, the number of vertices outside $X$ which are adjacent to $X$ is at least $c|X|$. 
Expander graphs share many of the properties of random regular graphs, and this is an important reason why they are of great interest in theoretical computer science. There are many excellent articles on expander graphs ranging from the very concise \cite{sarnak-expander} to the seriously comprehensive \cite{hoory-linial-widgerson}.

A key issue is that of constructing explicit expander graphs, and in particular that of constructing \emph{families} of expanders in which $k$ and $c$ are fixed but the number $n$ of vertices tends to infinity. Many constructions have been given, and several of them arise from Cayley graphs. Let $G$ be a finite group and suppose that $S = \{g_1^{\pm 1},\dots,g_k^{\pm 1}\}$ is a symmetric set of generators for $G$. The Cayley graph $\mathcal{C}(G,S)$ is the $2k$-regular graph on vertex set $G$ in which vertices $x$ and $y$ are joined if and only if $x y^{-1} \in S$. Such graphs provided some of the earliest examples of expanders \cite{lps,margulis-expander}. A natural way to obtain a family of such graphs is to take some large ``mother'' group $\tilde G$ admitting many homomorphisms $\pi$ from $\tilde G$ to finite groups, a set $\tilde S \subseteq \tilde G$, and then to consider the family of Cayley graphs $\mathcal{C}(\pi(\tilde G),\pi(\tilde S))$ as $\pi$ ranges over a family of homomorphisms. The work under discussion concerns the case $\tilde G = \SL_2(\Z)$, which of course admits homomorphisms $\pi_p : \SL_2(\Z) \rightarrow \SL_2(\F_p)$ for each prime $p$. For certain sets $\tilde S \subseteq \tilde G$, for example
\[ \tilde S = \{ \begin{pmatrix} 1 & 1 \\ 0 & 1\end{pmatrix}^{\pm 1}, \begin{pmatrix} 1 & 0 \\ 1 & 1\end{pmatrix}^{\pm 1}\}\] or
\[ \tilde S = \{ \begin{pmatrix} 1 & 2 \\ 0 & 1\end{pmatrix}^{\pm 1}, \begin{pmatrix} 1 & 0 \\ 2 & 1\end{pmatrix}^{\pm 1}\},\] 
spectral methods from the theory of automorphic forms may be used to show that $(\mathcal{C}(\pi_p(\tilde G),\pi_p(\tilde S)))_{\mbox{\scriptsize $p$ prime}}$ is a family of expanders. See \cite{lubotzky-survey} and the references therein. These methods depend on the fact that the group $\langle \tilde S \rangle$ has finite index in $\tilde G = \SL_2(\Z)$ and they fail when this is not the case, for example when \begin{equation}\label{123} \tilde S = \{ \begin{pmatrix} 1 & 3 \\ 0 & 1\end{pmatrix}^{\pm 1}, \begin{pmatrix} 1 & 0 \\ 3 & 1\end{pmatrix}^{\pm 1}\}.\end{equation} In \cite{lubotzky-survey} Lubotzky asked whether the corresponding Cayley graphs in this and other cases might nonetheless form a family of expanders, the particular case of \eqref{123} being known as his ``1-2-3 question''. The paper of Bourgain and Gamburd under discussion answers this quite comprehensively, showing that all that is required is that the group generated by $\tilde S$ is not \emph{virtually abelian} (contains a finite index abelian subgroup). We will sketch the proof in the case that $\tilde S$ generates a nonabelian free subgroup of $\SL_2(\Z)$. This is essentially the most general case, since the kernel of the natural homomorpism from $\langle \tilde S\rangle$ to $\SL_2(\F_2) \cong \mbox{Sym}(3)$ is free and has index at most $6$ in $\langle \tilde S\rangle$.

\begin{theorem}[Bourgain -- Gamburd]\label{bg-theorem}
Let $\tilde G = \SL_2(\Z)$ as above and suppose that $\tilde S$ is a finite symmetric set generating a free subgroup of $\SL_2(\Z)$. Then \[ (\mathcal{C}(\pi_p(\tilde G),\pi_p(\tilde S)))_{\mbox{\scriptsize $p$ \textup{prime}}}\] is a family of expanders.
\end{theorem}

The notation we have introduced here is rather cumbersome, so let us write $\Gamma_p := \mathcal{C}(\pi_p(\tilde G),\pi_p(\tilde S))$. For concreteness we will focus on the special case $\tilde S = \{A,A^{-1},B,B^{-1}\}$, where $A = \begin{pmatrix} 1 & 3 \\ 0 & 1\end{pmatrix}$ and $B = \begin{pmatrix} 1 & 0 \\ 3 & 1\end{pmatrix}$ are the matrices relevant to Lubotzky's 1-2-3 question. The argument is almost identical in any other case. In this case, then, $\Gamma_p$ is the graph on vertex set $\SL_2(\F_p)$ in which $x$ is joined to $y$ if and only if $xy^{-1}$ is one of the elements $A, A^{-1}, B$ or $B^{-1}$ considered modulo $p$.  Supposing that $p > 3$, each of these graphs is $4$-regular. The number of verices in $\Gamma_p$ is $n := |\SL_2(\F_p)| = p(p^2 - 1)$.

The reader may be interested to see a proof, using the ``ping-pong'' technique of Felix Klein, that that the subgroup of $\SL_2(\Z)$ generated by these $A$ and $B$ is indeed free. Consider the natural action of $A$ and $B$ on the projective plane $\P^1(\Q)$. Write 
\[ X := \{(\lambda : 1) \in \P^1(\Q) : |\lambda| < 1\}\] and
\[ Y := \{(1 : \lambda) \in \P^1(\Q) : |\lambda| < 1\},\] 
and observe that $X$ and $Y$ are disjoint and \emph{jouent au ping pong}, that is to say
\[ A^{n}(X) \subseteq Y\qquad \mbox{for all $n \in \Z \setminus \{0\}$}\]  and 
\[ B^{n}(Y) \subseteq X\qquad \mbox{for all $n \in \Z \setminus \{0\}$}.\] (The origin of the name should be clear -- the ``players'' $A$ and $B$ hit the domains $X$ and $Y$ back and forth -- as should the preference for the French term rather than the cumbersome ``play table tennis with one another''.) If the group generated by $A$ and $B$ is not free, then some nontrivial reduced word in $A$ and $B$ is equal to the identity, where ``reduced word'' means a finite word of the form $\dots A^{m_1}B^{n_1} \dots A^{m_k} B^{n_k} \dots$ with $m_1,n_1,\dots,m_k,n_k \neq 0$. The conjugate of such a word by an appropriate power of $A$ will still be the identity and will now have the form $w = A^{m_1} B^{n_1} \dots A^{m_k} B^{n_k} A^{m_{k+1}}$ with $m_i, n_j \neq 0$. However by repeated application of the ping-pong properties we see that $w(X) \subseteq Y$, certainly an impossibility since $X$ and $Y$ are disjoint and $w$ is supposed to be the identity.

Following that slight digression let us focus once again on the Cayley graphs $\Gamma_p$, our aim being to prove that they form a family of expanders as $p$ ranges over the primes. To do this we begin by giving a spectral interpretation of the expansion property which we defined combinatorially above. For each $p$ we may consider the \emph{Laplacian} of the corresponding Cayley graph, that is to say the operator 
\[ \Delta : L^2(\SL_2(\F_p)) \rightarrow L^2(\SL_2(\F_p))\] defined by
\[ \Delta f(x) := f(x) - \textstyle\frac{1}{4}\displaystyle( f(A x) + f(A^{-1}x ) + f(Bx) + f(B^{-1} x)).\]
The eigenvalues of the Laplacian lie in the interval $[0,2]$. Zero is certainly an eigenvalue, since $\Delta 1 = 0$. Write the eigenvalues in ascending order as $0 = \lambda_0 \leq \lambda_1 \leq \dots \leq \lambda_{n-1}$. It turns out the expansion properties of the graph $\Gamma_p$ (in fact of \emph{any} regular graph) are intimately connected with the size of the second-smallest eigenvalue $\lambda_1 = \lambda_1(\Gamma_p)$. The precise relation between the combinatorial property of expansion and this spectral property is discussed in Section 2 of \cite{hoory-linial-widgerson}, but for our purposes we need only remark that it suffices to show that the second-smallest eigenvalue $\lambda_1(\Gamma_p)$ is bounded away from zero independently of $n$ (in fact, this is also a necessary condition for expansion). The term \emph{spectral gap} is used to describe this property: there is a gap at the bottom of the spectrum in which there are no eigenvalues apart from zero.

To try to show that there is a spectral gap, consider the operator 
\[ T : L^2(\SL_2(\F_p)) \rightarrow L^2(\SL_2(\F_p))\] given by $T := 4 (\id - \Delta)$, that is to say
\[ Tf(x) := f(A x) + f(A^{-1}x ) + f(Bx) + f(B^{-1} x).\]
The matrix\footnote{With respect to the basis of $\SL_2(\F_p)$ consisting of the functions $\mathbf{1}_t : \SL_2(\F_p) \rightarrow \C$ defined by $\mathbf{1}_t(x) = 1$ if $x = t$ and $0$ otherwise.} of $T$ is same thing as the adjacency matrix of the graph $\Gamma_p$, that is to say the matrix whose $xy$ entry is $1$ if $x \sim y$ and zero otherwise. The eigenvalues of $T$ are of course $\mu_i = 4(1 - \lambda_i)$, $i = 0,\dots,n-1$, and it is a very well-known and easy to establish fact that the $2m$th moment $\sum_{i=0}^{n-1} \mu_i^{2m}$ is equal to $n$ times $W_{2m}$, the number of closed walks of length $2m$ from the identity to itself. It follows that we have
\begin{equation}\label{trace-formula} W_{2m} = \frac{1}{n}4^{2m} \left( 1 + \sum_{i = 1}^{n-1} (1 - \lambda_i)^{2m}\right).\end{equation}
Note in particular that $W_{2m} \geq \frac{1}{n}4^{2m}$, since all the terms are non-negative. At first glance it looks as though the only way to use \eqref{trace-formula} to bound $\lambda_1$ away from zero would be to get rather precise estimates on $W_{2m}$, and in particular one would at the very least want to show that $W_{2m} < \frac{2}{n}4^{2m}$. However a remarkable observation, used earlier in related contexts by Sarnak and Xue \cite{sarnak-xue} and Gamburd \cite{gamburd-56}, comes into play. This is that any eigenspace of the Laplacian is $\SL_2(\F_p)$-invariant, where the action of $\SL_2(\F_p)$ on $L^2(\SL_2(\F_p))$ is the right-regular one given by $g \circ f(x) := f(xg)$. In other words, any such eigenspace has the structure of a representation of $\SL_2(\F_p)$ and thus, by basic representation theory, decomposes into irreducible representation of $\SL_2(\F_p)$. But by a classical theorem of Frobenius all such representations have dimension at least $(p-1)/2 \sim n^{1/3}$. This means that $\lambda_1 = \lambda_2 = \dots = \lambda_l$ for some $l \sim n^{1/3}$, and hence from \eqref{trace-formula} we in fact have the bound
\begin{equation}\label{multiplicity-bound} W_{2m} \gg \frac{1}{n^{2/3}}4^{2m}(1 - \lambda_1)^{2m}.\end{equation}
This enables a meaningful spectral gap (lower bound on $\lambda_1$) to be obtained from somewhat weaker upper bounds on $W_{2m}$.

The main new content of \cite{bourgain-gamburd-sl2}, then, is to obtain those upper bounds on $W_{2m}$, the number of walks of length $2m$ starting and finishing at the identity, for appropriate $m$. A nice way of thinking about these walks is in terms of \emph{convolution powers} of the probability measure
\[ \nu := \textstyle\frac{1}{4}\displaystyle( \delta_{A} + \delta_{A^{-1}} + \delta_{B} + \delta_{B^{-1}} )\] on $\SL_2(\F_p)$, where $\delta_g(x) = n$ if $x = g$ and $0$ otherwise. This measure $\nu$ is a very singular or ``spiky'' probability measure, supported on just the four points $A, A^{-1}, B$ and $B^{-1}$. Now the convolution
\[ \nu^{(2)} := \nu \ast \nu(x) := \E_{y \in \SL_2(\F_p)} \nu(xy^{-1})\nu(y)\] is supported on words of length at most two in $A, A^{-1}, B$ and $B^{-1}$, or alternatively those $x$ in the graph $\Gamma_p$ which can be reached from the identity by a path of length two, the value of $\nu \ast \nu(x)$ being $4^{-2}n$ times the number of paths of length two from the identity to $x$. Similarly higher convolution powers $\nu^{(m)}(x) := \nu \ast \dots \ast \nu(x)$ give $4^{-m}n$ times the number of paths of length $m$ from the identity to $x$. The idea of the proof is to examine these convolution powers, showing that they become progressively more ``spread out'' until, for suitable $m$, $\nu^{(2m)}$ vaguely resembles the uniform measure $\mathbf{1}$ which assigns weight one to each point of $\SL_2(\F_p)$. Then, in particular, $\nu^{(2m)}(0) \sim 1$, meaning that $W_{2m} \sim 4^{2m}/n$. Combined with \eqref{multiplicity-bound}, this is enough to establish the desired spectral gap.

The notion of a probability measure $\mu$ on $\SL_2(\F_p)$ being ``spread out'' may be quantified using the $L^2$-norm \[ \Vert \mu \Vert_2 := \big(\E_{x \in \SL_2(\F_p)} \mu(x)^2\big)^{1/2}.\] The $L^2$-norm of a delta measure $\delta_g$ is $n^{1/2}$, which is huge, whilst that of the uniform measure $\mathbf{1}$ is equal to one, the smallest value possible by the Cauchy-Schwarz inequality. It is not hard to show that convolution cannot increase the $L^2$-norm, and so we have the chain of inequalities
\begin{equation}\label{sequence} n^{1/2} = \Vert \nu^{(1)}\Vert_2 \geq \Vert \nu^{(2)}\Vert_2 \geq \dots.\end{equation}
The aim is to show that this sequence is, in fact, rather rapidly decreasing. Roughly speaking one shows that 
\begin{equation}\label{spread-out-final} \Vert \nu^{(m_1)} \Vert_2 \approx 1\end{equation} for some $m_1 \approx C_1\log p$; this $m_1$ turns out to be an appropriate choice to substitute into \eqref{multiplicity-bound} in order to reach the desired conclusions.

It turns out that this sequence gets off to a rather good start. This is a consequence of an observation of Margulis \cite{margulis-combinatorica}, namely that the freeness of the subgroup of $\SL_2(\Z)$ generated by $A$ and $B$ persists to some extent even when reduced modulo $p$. Indeed let us take a reduced word $w= A^{m_1} B^{n_1} \dots A^{m_k} B^{n_k}$ with $m_1,\dots,m_k, n_1,\dots,n_k \neq 0$ and suppose that this equals the identity when reduced modulo $p$, that is to say in $\SL_2(\F_p)$. Lifting back up to $\SL_2(\Z)$ we have
\[ \tilde w = A^{m_1} B^{n_1} \dots A^{m_k} B^{n_k} \equiv \id \pmod{p}.\]
But the freeness of the lifted group means that $\tilde w \neq \id$, and thus in order to be congruent to the identity mod $p$ the matrix $\tilde w$ must have at least one entry of size at least $p-1$. But by some simple matrix inequalities this is impossible provided that 
\[ |m_1| + |n_1| + \dots + |m_k | + |n_k| < c \log p\] for some absolute constant $c > 0$.

It follows that the subgroup of $\SL_2(\F_p)$ generated by $A$ and $B$ is ``free up to words of length $c \log p$''. In terms of the Cayley graphs $\Gamma_p$ this means that up to retracing steps there is a unique walk of length $m$ between the identity and $x$ for any $x \in \SL_2(\F_p)$ and for any $m \leq m_0 := c\log p/2$. This implies that the measures $\nu^{(m)}$, $m \geq m_0$ are already rather spread out. To quantify this (and in particular to deal with the issue of ``retracing steps'') a result of Kesten concerning random walks in the free group may be applied. The conclusion is that\begin{equation}\label{initial-spread}\Vert \nu^{(m_0)} \Vert_2 \ll n^{1/2 - \gamma}\end{equation} for some $\gamma > 0$. This is good progress on the way to \eqref{spread-out-final} and represents a significant improvement on the initial bound $\Vert \nu^{(1)} \Vert_2 = n^{1/2}$.

It is convenient to imagine, for the rest of the argument, that all probability measures $\mu$ on $G$ have the form $\mu(x) = \frac{n}{|A|} 1_A(x)$ for some set $A \subseteq G$, the ``support'' of $\mu$. Whilst this is clearly not true, various (somewhat technical) decompositions into level sets may be used to reduce to this case. For such a measure we have \[ \Vert \mu \Vert_2 = (n/|A|)^{1/2},\] and so the bound \eqref{initial-spread} corresponds to $|A| \gg n^{2\gamma}$, certainly a reasonable level of spreadoutness.

The rest of the argument, which constitutes the heart of the paper, involves examining the convolution powers between $\nu^{(m_0)}$ and $\nu^{(m_1)}$ for a suitable $m_1 \sim C_1 \log p$, the aim being to establish \eqref{spread-out-final}. An application of the ``dyadic pigeonholing argument'', used to great effect by Bourgain in many papers, is employed: if $\Vert \nu^{(m_1)}\Vert_2$ is much larger than 1, this means that the sequence \eqref{sequence} cannot decay too rapidly between $\nu^{(m_0)}$ and $\nu^{(m_1)}$ and so there must be two convolution powers $\nu^{(m)}$ and $\nu^{(2m)}$, $m_0 \leq m < m_1$, such that $\Vert \nu^{(2m)} \Vert_2 \approx \Vert \nu^{(m)} \Vert_2$. Let us be deliberately vague about the meaning of $\approx$ here.

Suppose that $\nu^{(m)}(x) = \frac{n}{|A|}1_A(x)$ for some set $A \subseteq G$. Noting that $\nu^{(2m)} = \nu^{(m)} \ast \nu^{(m)}$, it is not hard to compute that the ratio \[ \Vert \nu^{(2m)}\Vert_2^2/\Vert \nu^{(m)} \Vert_2^2\] is actually equal to $|A|^{-3}$ times the number of quadruples $a_1,a_2,a_3,a_4 \in A^4$ with $a_1 a_2 = a_3 a_4$. This may be compared with condition (4) in the list of properties which are known to roughly characterise approximate groups. Thus, being even rougher at this point, \begin{equation}\label{nu-h}\nu^{(m)} \sim \frac{1}{H}1_H\end{equation} for some approximate group $H \subseteq \SL_2(\F_p)$. Note that the rough equivalence of (4) and other, more flexible definitions such as Definition \ref{approx-gp-def} is one of the deeper equivalences mentioned in \S \ref{approx-groups-sec}, being reliant on the nonabelian Balog-Szemer\'edi-Gowers theorem of Tao \cite{tao-noncommutative}.

If $H$ is already all of $\SL_2(\F_p)$ then \eqref{nu-h} is telling us that $\nu^{(m)}$ is close to the uniform distribution, in which case so is $\nu^{(m_1)}$, hence \eqref{spread-out-final} is established and we are done. If not then we apply Helfgott's result, Theorem \ref{helfgott-theorem}, to conclude that $H$ is essentially upper-triangular, and hence that $\nu^{(m_0)}$ has significant mass on an upper-triangular subgroup of $\SL_2(\F_p)$. 

The support of $\nu^{(m_0)}$, however, consists of words of length at most $m_0$ in the generators $A, A^{-1}, B$ and $B^{-1}$ and, as we stated, these elements behave freely up to words of this length. This is highly incompatible with upper-triangularity, which in particular implies that we always have the commutator relation\footnote{In other words, upper-triangular subgroups of $\SL_2(\F_p)$ are 2-step solvable.} \begin{equation}\label{commutator}[[g_1,g_2],[g_3,g_4]] = \id.\end{equation} A pleasant group-theoretic argument formalises this incompatibility and allows one to show that any set of words of length at most $m_0$ in the generators $A, A^{-1}, B$ and $B^{-1}$ satisfying \eqref{commutator} has size at most $m_0^6$. This represents a tiny proportion of the set of all such words, which (counted with multiplicity at least) has cardinality $4^{m_0}$. This contradiction finishes the sketch proof of Theorem \ref{bg-theorem}.\hfill $\Box$

Before moving on, we wish to record, for use in the next section, a further observation concerning the measures $\nu^{(m)}$. We sketched a proof that $\Vert \nu^{(m_1)} \Vert_2 \approx 1$ for some $m_1 \sim C_1\log p$, that is to say $\nu^{(m_1)}$ vaguely resembles the uniform distribution on $\SL_2(\F_p)$. By taking further convolutions and using the fact that irreducible representations have large degree once more, this may be bootstrapped to show that $\nu^{(m)}$ becomes exponentially well uniformly-distributed:
\begin{equation}\label{conv} \nu^{(m)}(x) = 1 + O(ne^{-cm})\end{equation} for some absolute $c > 0$ and for all $m$. Alternatively, such a statement can be deduced directly from the spectral gap property, as is done for example in \cite[\S 3.3]{bourgain-gamburd-sarnak}.

It is interesting to ask whether the arguments might adapt to deal with Cayley graphs on $\SL_n(\F_p)$ with $n \geq 3$. A recent paper of Bourgain and Gamburd \cite{bourgain-gamburd-pn-ii} shows that this is the case when $n = 3$. The argument is, in large part, quite similar to the above, except of course that Helfgott's theorem on approximate subgroups of $\SL_2(\F_p)$ must be replaced by his more difficult result \cite{helfgott-2} on approximate subgroups of $\SL_3(\F_p)$. There is one significant extra difficulty, however, which is that there are proper subgroups of $\SL_3(\F_p)$ which are not close to upper-triangular, an obvious example being a copy of $\SL_2(\F_p)$. To deal with this a deep algebro-geometric result of Nori \cite{nori} is brought into play, which states that any proper subgroup of $\SL_3(\F_p)$, $p$ sufficiently large, must satisfy a non-trivial polynomial equation. To obtain a contradiction, it must be shown that the set of words of length $m_0$ in the generators $A$ and $B$ (say) does not concentrate on the corresponding subvariety of $\SL_3(\C)$, and here techniques from the theory of random matrix products and a certain amount of ``quantitative algebraic geometry'' are brought into play.

\section{Nonlinear sieving problems}

In this section we discuss work of Bourgain, Gamburd and Sarnak \cite{bourgain-gamburd-sarnak}. The goal of \emph{sieve theory}, traditionally viewed as a part of analytic number theory, is to find prime numbers or at least to say something about them. Historically, the sieve arose through work of Brun and Merlin on the twin prime problem, that is to say the problem of finding infinitely may primes $p$ such that $p + 2$ is also prime. Whilst this remains a famous open problem, approximations to it have been found. For example, Brun established the following result.

\begin{theorem}[Brun]
There are infinitely many integers $n$ such that $n(n+2)$ has at most $9$ prime factors. 
\end{theorem}

Much later, Chen \cite{chen} replaced $9$ by $3$. One way of stating this type of result is as follows: there are infinitely many $n$ for which both $n(n+2)$ is a 3-\emph{almost prime}, that is to say a positive integer with at most 3 prime factors. 

The aim of \cite{bourgain-gamburd-sarnak} is to discover almost primes in more exotic locales, and specifically in \emph{orbits} of linear groups. We will sketch a proof of the following result.

\begin{theorem}[Bourgain-Gamburd-Sarnak]\label{bgs-special-case}
Let $A$ and $B$ be two matrices in $\SL_2(\Z)$ generating a free subgroup. Then there is some $r$ such that this group contains infinitely many $r$-almost prime matrices \textup{(}matrices, the product of whose entries is $r$-almost prime\textup{)}.
\end{theorem}

Henceforth we shall say ``almost prime'' instead of ``$r$-almost prime for some $r$''. We remark that in the specific case we focussed on in the last section, when $A = \begin{pmatrix} 1 & 3 \\ 0 & 1\end{pmatrix}$ and $B = \begin{pmatrix} 1 & 0 \\ 3 & 1\end{pmatrix}$, the theorem as stated follows from classical sieve theory of the type used to prove Brun's theorem. Indeed (for example) $A^nBA = \begin{pmatrix} 9n+1 & 30n+3 \\ 3 & 10\end{pmatrix}$, and the product of the entries here is $2 \cdot 3^2 \cdot 5 \cdot (9n+1) \cdot (10n+1)$, which will be almost prime for infinitely many $n$ by a simple variant of Brun's analysis. The issue here is that the subgroup generated by $A$ and $B$ contains unipotent elements (in this case both $A$ and $B$ are themselves unipotent). 

We start with a (very) elementary discussion of what a sieve is. Suppose one has a finite set $X$ of integers and that one wishes to find primes or almost primes in $X$. The most na\"{\i}ve way to do this would be to try to adapt the sieve of Eratosthenes, using the inclusion-exclusion principle to compute
\[ \#\{ \mbox{primes in $X$}\} = |X| - |X_2| - |X_3| - |X_5| - \dots + |X_6| + |X_{10}| + |X_{15}| + \dots - |X_{30}| - \dots\] where $X_q$ is the set of elements of $X$ which are divisible by $q$. Unfortunately it is well-known that, even when $X$ is an extremely simple set such as $\{1,\dots,n\}$, it is not generally possible to evaluate $|X_q|$ sufficiently accurately to avoid the error terms in this long sum blowing up. In this simple case just mentioned, for example, we have $|X_q| = \lfloor n/q\rfloor$. However the floor function is rather unpleasant and it is tempting to write instead $|X_q| = n/q + O(1)$, but then one finds that there are so many $O(1)$ errors that the sieve of Eratosthenes becomes useless.

By and large, \emph{sieve theory} is concerned with what it is possible to say about primes or almost primes in $X$ given ``reasonably nice'' information about the size of the sets $X_q$. Although the sieve of Eratosthenes is bad, other sieves fare rather better. These other sieves are generally cleverly weighted versions of the sieve of Eratosthenes, but we will not dwell upon their construction here. A typical example of ``reasonably nice'' information about $|X_q|$ would be
\[ |X_q| = \beta(q)|X| + r_q\] for all squarefree $q \leq |X|^{\gamma}$, where $\beta(q)$ is some pleasant multiplicative function and the error $r_q$ is small in the sense that $|r_q| \ll |X|^{1 - \gamma}$ for some $\gamma > 0$.  For example, if $X = \{1,\dots,n\}$ then this is true with $\beta(q) = 1/q$ and for any $\gamma \leq 1$. 

The fundamental theorem of the combinatorial sieve states, roughly speaking, that such information is enough to find almost primes in $X$; in fact, one can even estimate the number of almost primes. What is meant by ``almost prime'' -- that is, how many prime factors these numbers will have -- depends on how large we can take $\gamma$ as well as on the so-called dimension of the sieve, which has to do with the average size of the quantities $\beta$. We will not delay ourselves by expanding upon the details here. Let us instead refer the reader to \cite{bourgain-gamburd-sarnak} for the precise formulation convenient to the application there and to the book \cite{iwaniec-kowalski} or the unpublished notes \cite{iwaniec} for a more wide-ranging discussion of sieves in general with full proofs.

All we shall take from the preceding discussion is the notion that, given a finite set $X$ to be sieved in order to locate almost primes, we should be looking for good asymptotics for the size of the sets $|X_q|$, $q$ squarefree. Returning to Theorem \ref{bgs-special-case}, the first obvious question to answer is that of what the set $X$ to be sieved should be. The set in which we wish to find almost primes is
\[ \mathcal{A} := \{ x_1 x_2 x_3 x_4 :  \begin{pmatrix} x_1 & x_2 \\  x_3& x_4\end{pmatrix}\in \langle A, B\rangle\}.\]  Now $\mathcal{A}$ is of course an infinite set of integers. Rather than truncate in the usual way and take $X = \mathcal{A} \cap \{1,\dots,N\}$, it is much more natural to truncate in a manner that respects the group structure more. This we do by taking
\[ X:=  \{ x_1 x_2 x_3 x_4 :  \begin{pmatrix} x_1 & x_2 \\  x_3& x_4\end{pmatrix}\in \Sigma_m(A,B)\},\] where
\[  \Sigma_m(A,B) = \{ U_1U_2 \dots U_m : U_i \in \{A,A^{-1},B,B^{-1}\} \}\] is the set of words of length $m$ in $A, A^{-1}, B$ and $B^{-1}$ and $X$ is counted with multiplicity so that $|X| = 4^m$.

Suppose that $p$ is a prime. Then $|X_p|$ is equal to the number of words $w \in \Sigma_m(A,B)$, counted with multiplicity, which, when reduced modulo $p$, give rise to a matrix in $\SL_2(\F_p)$ with at least one zero entry.  Writing $S \subseteq \SL_2(\F_p)$ for the set of such matrices, it is easy to compute that $|S| = 2(2p-1)(p-1)$. Now the number of words $w \in \Sigma_m(A,B)$ which reduce modulo $p$ to some $x \in \SL_2(\F_p)$ is, in the notation of the last section, precisely $\frac{1}{n}|X| \nu^{(m)}(x)$, and so
\[ |X_p| = \frac{1}{n} |X| \sum_{x \in S} \nu ^{(m)}(x).\]
However at the end of the last section we saw\footnote{Either as a byproduct of the proof, or a consequence, of the expansion property of the family of Cayley graphs $\Gamma_p = \mathcal{C}(\pi_p(\mathcal{G}), \pi_p(\tilde S))$.} that $\nu^{(m)}(x)$ becomes very close to $1$. In fact, in \eqref{conv} we noted the bound $\nu^{(m)}(x) = 1 + O(n e^{-cm})$. Using this we obtain
\[ |X_p| = \beta(p)|X| + r_p \] where $\beta(p) := 2(2p-1)/p(p+1)$ and $|r_p| = |X|^{1 - \gamma}$ for some $\gamma > 0$. 

Thus the expansion property of the Cayley graphs $(\mathcal{C}(\pi_p(\tilde G), \pi_p(\tilde S)))_{\mbox{\scriptsize $p$ prime}}$ gives exactly the kind of information that can be input into the combinatorial sieve!

There is, however, a very major caveat. What we have just said applies only to $X_p$ when $p$ is a prime, and for the sieve one must understand $X_q$ when $q$ is a general squarefree number. To do this requires the establishment of Theorem \ref{bg-theorem} for the family $(\mathcal{C}(\pi_q(\tilde G), \pi_q (\tilde S)))_{q}$, where now $q$ ranges over all squarefrees and not just over primes. The broad scheme of the proof is the same, but every single ingredient must be generalised to the more general setting, starting from the classification of approximate subrings of $\Z/q\Z$. The situation here is more complicated because this ring will in general have many approximate subrings, namely $\Z/q'\Z$ with $q' | q$. One of the main technical results of \cite{bourgain-gamburd-sarnak} (occupying some 20 pages) is the statement that, very roughly speaking, these are the only approximate subrings of $\Z/q\Z$. Although this is a deeply technical argument of a type that this author would struggle to summarise meaningfully even to an expert audience, it might be compared with the 92-page proof \cite{bourgain-zq} of the corresponding assertion without the squarefree assumption on $q$. Thankfully\footnote{This is one of the most extraordinarily long and technical arguments the author has ever seen. The theory of approximate rings when there are many zero-divisors seems to be very difficult.} this is not required for the present application. Once the classification of approximate subrings of $\Z/q\Z$ for $q$ squarefree is in place a suitable analogue of Helfgott's argument is applied to roughly classify approximate subgroups of $\SL_2(\Z/q\Z)$. Even the statement of this result (Proposition 4.3 in the paper) is rather technical. Finally, the majority of the argument outlined in the last section in the case $q$ prime goes over without substantial change.

This concludes our discussion of the proof of Theorem \ref{bgs-special-case}. To conclude this survey, we wish to mention a beautiful application, mentioned in the original paper \cite{bourgain-gamburd-sarnak} and in other articles such as \cite{sarnak-survey}, of these nonlinear sieving ideas. This has to do with Apollonian packings such as the one in the attractive image above.
\begin{figure}
\includegraphics{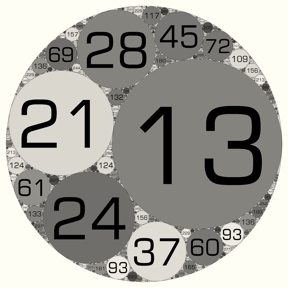}
\caption{Apollonian circle packing}
\label{figure-1}
\end{figure}

For a very pleasant and gentle introduction to Apollonian packings, see \cite{austin}. Referring to Figure \ref{figure-1}, inside each circle is an integer which represents the curvature of that circle, or in other words the reciprocal of the radius. Some of the number theory associated with the integers that arise in this way is discussed in the letter \cite{sarnak-letter} where, for example, it is shown that infinitely many of these curvatures are prime and in fact that there are infinitely many touching pairs of circles with prime curvature. 

Now a pleasant exercise in Euclidean geometry gives a theorem of Descartes, namely that the relation between the four integers $a_1,a_2,a_3,a_4$ inside four mutually touching circles is given by 
\begin{equation}\label{quadratic} 2(a_1^2 + a_2^2 + a_3^2 + a_4^2) = (a_1 + a_2 + a_3 + a_4)^2.\end{equation} Examples of quadruples $(a_1,a_2,a_3,a_4)$ which are related in this way and easily visible in the picture are $(13,21,24,124)$ and $(13,24,37,156)$.

Take a quadruple $(C_1,C_2,C_3,C_4)$ of touching circles with curvatures \[ (a_1,a_2,a_3,a_4) = (13,21,24,124).\] There is another circle $C'_1$ tangent to $C_2,C_3$ and $C_4$, and it has curvature $a'_1 = 325$. To find a general relation between $a_1$ and $a'_1$ we may note that $a_1,a'_1$ are roots of \eqref{quadratic} regarded as a quadratic in $a_1$ and thereby obtain the relation
\[ a'_1 = -a_1 + 2a_2 + 2a_3 + 2a_4.\] This may of course be written as
\[ \begin{pmatrix} a'_1 \\ a_2 \\ a_3 \\ a_4\end{pmatrix} = \begin{pmatrix}  -1 & 2 & 2 & 2 \\ 0 & 1 & 0 & 0 \\ 0 & 0 & 1 & 0 \\ 0 & 0 & 0 & 1  \end{pmatrix}\begin{pmatrix} a_1 \\ a_2 \\ a_3 \\ a_4\end{pmatrix}.\]
That is, if one starts with some fixed vector such as $x_0 = (13,21,24,124)$ then one may obtain another quadruple of curvatures of circles in the Apollonian packing by applying the matrix 
\[ S_1 :=  \begin{pmatrix}  -1 & 2 & 2 & 2 \\ 0 & 1 & 0 & 0 \\ 0 & 0 & 1 & 0 \\ 0 & 0 & 0 & 1  \end{pmatrix}.\]
By playing the same game with $C'_2,C'_3$ and $C'_4$ we can make the same assertion with the matrices
\[ S_2 :=  \begin{pmatrix}  1 & 0 & 0 & 0 \\ 2 & - 1 & 2 & 2 \\ 0 & 0 & 1 & 0 \\ 0 & 0 & 0 & 1 \end{pmatrix},\]
\[ S_3 :=  \begin{pmatrix}  1 & 0 & 0 & 0 \\ 0 & 1 & 0 & 0 \\ 2 & 2 & - 1& 2 \\ 0 & 0 & 0 & 1 \end{pmatrix}\]
and
\[ S_4 :=  \begin{pmatrix}  1 & 0 & 0 & 0 \\ 0 & 1 & 0 & 0 \\ 0 & 0 & 1 & 0 \\ 2 & 2 & 2 & -1 \end{pmatrix}.\]
This leads naturally to consideration of the orbit $\langle \tilde S \rangle x_0 \subseteq \Z^4$, where 
\[ \tilde S := \{ S_1, S_2, S_3, S_4\},\] every vector in which consists entirely of curvatures of circles in the Apollonian packing. This puts us in a situation very similar to that studied in Theorem \ref{bgs-special-case}, except now we appear to be dealing with a subgroup of $\GL_4(\Z)$ rather than of $\SL_2(\Z)$.

It turns out, however, that this situation is essentially a two-dimensional one in disguise, and for this we need to add to the list of areas of mathematics we touch upon by hinting at Lie theory and special relativity! The matrices $S_1,S_2,S_3,S_4$ belong to $\SO_F(\Z)$, the subgroup of $\GL_4(\Z)$ consisting of $4 \times 4$ matrices with determinant one which preserve the quadratic form $F(\vec{x}) = 2(x_1^2 + x_2^2 + x_3^2 + x_4^2) - (x_1 + x_2 + x_3 + x_4)^2$ (cf. \eqref{quadratic}). By the standard theory of quadratic forms (over $\R$) this is equivalent to the Lorentz form $L(\vec{y}) = y_1^2 + y_2^2 + y_3^2 - y_4^2$, and so we may identify $\SO_F(\R)$ with the orthogonal group $\SO(3,1)$ preserving this latter form. But it is very well-known that this group admits $\mbox{SU}(2)$ as a double cover: this is because the set $\{\vec{y} : L(\vec{y}) = -1\}$ may be identified with the set of $2 \times 2$ hermitian matrices $M$ with determinant 1 via
\[ (y_1,y_2,y_3,y_4) \mapsto \begin{pmatrix} y_4 + y_3 & y_1 - i y_2 \\ y_1 + iy_2 & y_4 - y_3\end{pmatrix},\] and so any $P \in \mbox{SU(2)}$ gives rise to an element of $\SO(3,1)$ via the transformation $M \mapsto P M P^*$.

By lifting to this double cover the group $\langle \tilde S\rangle$ can be lifted to a subgroup of $\SL_2(\Z[i])$. The proof of Theorem \ref{bgs-special-case} goes through with relatively minimal changes, although once again the group generated by $S_1,S_2,S_3$ and $S_4$ contains unipotents and so, if the aim is simply to find infinitely many circles or pairs/quadruples of touching circles with almost-prime curvatures, more elementary approaches work just as well. 
Those elementary approaches do not, however, give sharp quantitative results, whereas the techniques we have sketched do. To explain one such result, imagine Figure \ref{figure-1} being generated as follows. Start with the outer circle (which has curvature $-6$) and the three largest inner circles, with curvatures 13, 21 and 24. This is the \emph{first generation}. The second generation consists of those circles touching three from the first generation: they have curvatures 28,37,61 and 124. The third generation contains those new circles touching three circles from either the first or the second generations: these have curvatures 45, 60, 69, 93, 124, 132, 133, 156,  220, 292, 301 and 325. Carry on in this vein: the $n$th generation will contain $4 \cdot 3^{n-2}$ circles.

\begin{theorem}[Bourgain, Gamburd, Sarnak]
The number of circles at generation $n$ which have prime curvature is bounded by $C 3^n/n$, for some absolute constant $C$.
\end{theorem}

We conclude by remarking that there are some very interesting unsolved questions connected with Apollonian packings \cite{lagarias}. In that paper the very interesting question is raised of whether, in Figure \ref{figure-1}, a positive proportion of all positive integers appear as curvatures. J.~Bourgain has recently indicated to me that he and Elena Fuchs have obtained new information on this question. See also \cite{kont-oh} for an asymptotic formula for the number of circles in the packing of curvature at most $X$. It seems that the question of describing this set of integers more precisely remains open: are they given, from some point on, by finitely many congruence conditions?

\section{Acknowledgements}

I am very grateful to Cliff Reiter for his permission to use the image of an Apollonian packing. In addition to the cited references, the notes from the fourth Marker Lecture of Terence Tao were very useful in preparing what is written here. I thank Jean Bourgain, Emmanuel Breuillard, Alex Gamburd, Harald Helfgott, Alex Kontorovich, Peter Sarnak and Terry Tao for their comments on an earlier draft of this survey. Finally, I am very grateful to Lilian Matthiesen and Vicky Neale for correcting a number of typographical errors.


\bibliographystyle{amsplain}

\begin{thebibliography}{99}

\bibitem{austin} D.~Austin, \emph{When kissing involves trigonometry,} AMS monthly feature article, available at \texttt{http://www.ams.org/featurecolumn/archive/kissing.html}.

\bibitem{bourgain-zq} J.~Bourgain, \emph{The sum-product theorem in $\Z_q$ with $q$ arbitrary,} J. Anal. Math. \textbf{106} (2008), 1--93.

\bibitem{bourgain-gamburd-sl2} J.~Bourgain and A.~Gamburd, \emph{Uniform expansion bounds for Cayley graphs of $\SL_2(\F_p)$,} Ann. Math. \textbf{167} (2008), 625--642.

\bibitem{bourgain-gamburd-pn} \bysame, \emph{Expansion and random walks in $\SL_d(\Z/p^n\Z)$. I,} J. Eur. Math. Soc. (JEMS) \textbf{10} (2008), no. 4, 987--1011. 

\bibitem{bourgain-gamburd-pn-ii} \bysame \emph{Expansion and random walks in $\SL_d(\Z/p^n\Z)$, II,} with an appendix by J.~Bourgain, J. Eur. Math. Soc. (JEMS) \textbf{11} (2009), 1057--1103.

\bibitem{bourgain-gamburd-sarnak} J.~Bourgain, A.~Gamburd and P.~Sarnak, \emph{Affine linear sieve, expanders, and sum-product,} preprint.

\bibitem{bourgain-glibichuk-konyagin}  J.~Bourgain, A.~A.~Glibichuk and S.~V.~Konyagin, \emph{Estimates for the number of sums and products and for exponential sums in fields of prime order,} J. London Math. Soc. (2) \textbf{73} (2006), no. 2, 380--398. 

\bibitem{bkt} J.~Bourgain, N.~H.~Katz and T.~C.~Tao, \emph{A sum-product estimate in finite fields and applications,} Geom. Funct. Anal. (GAFA) \textbf{14} (2004), 27--57.

\bibitem{bourgain-konyagin} J.~Bourgain and S.~Konyagin, \emph{Estimates for the number of sums and products and for exponential sums over subgroups in fields of prime order,} C.~R.~Acad. Sci. Paris, Ser. I, \textbf{337} (2003), 75--80.

\bibitem{bg1} E.~Breuillard and B.~J.~Green, \emph{Approximate groups, I: the torsion-free nilpotent case,} preprint available at \texttt{http://arxiv.org/abs/0906.3598}.

\bibitem{bg2} \bysame, \emph{Approximate groups, II: the solvable linear case,} preprint available at \texttt{http://arxiv.org/abs/0907.0927}.

\bibitem{bg3} \bysame, \emph{Approximate groups, III: the bounded linear case,} in preparation.

\bibitem{chang-sl3} M.-C.~Chang, \emph{Product theorems in $\SL_2$ and $\SL_3$,} J. Inst. Math. Jussieu \textbf{7} (2008), no. 1, 1--25. 

\bibitem{chen} J.~R.~Chen, \emph{On the representation of a larger even integer as the sum of a prime and the product of at most two primes,} Sci. Sinica \textbf{16} (1973), 157--176. 

\bibitem{edgar-miller} G.~A.~Edgar and C.~Miller, \emph{Borel subsets of the reals,} Proc. Amer. Math. Soc. \textbf{131} (2003), 1121--1129.

\bibitem{elekes-kiraly} G.~Elekes and Z.~Kir\'aly, \emph{On the combinatorics of projective mappings,} J. Algebraic Combin. \textbf{14} (2001), no. 3, 183--197. 

\bibitem{erdos-szemeredi} P.~Erd\H{o}s and E.~Szemer\'edi, \emph{On sums and products of integers,} Studies in pure mathematics, 213--218, Birkh\"auser, Basel, 1983. 

\bibitem{eskin-mozes-oh} A.~Eskin, S.~Mozes and H.~Oh, \emph{On uniform exponential growth for linear groups,} Invent. Math. \textbf{160} (2005), no. 1, 1--30.

\bibitem{fisher-katz-tao} D.~Fisher, N.~H.~Katz and I.~Peng, \emph{On Freiman's theorem in nilpotent groups,} preprint available at \texttt{http://arxiv.org/abs/0901.1409}.

\bibitem{freiman-book} G.~A.~Fre\u{\i}man, \emph{Foundations of a structural theory of set addition,} Translations of Mathematical Monographs, \textbf{37}. American Mathematical Society, Providence, R. I., 1973. vii+108 pp.

\bibitem{gamburd-56} A.~Gamburd, \emph{On the spectral gap for infinite index ``congruence'' subgroups of $\SL_2(\Z)$}, Israel J. Math. \textbf{127} (2002), 157--200.

\bibitem{gowers1} W.~T.~Gowers, \emph{A new proof of Szemer\'edi's theorem for arithmetic progressions of length four,} Geom. Funct. Anal. \textbf{8} (1998), no. 3, 529--551. 

\bibitem{gowers-gafa} \bysame, \emph{Rough structure and classification,} GAFA 2000 (Tel Aviv, 1999). Geom. Funct. Anal. 2000, Special Volume, Part I, 79--117. 

\bibitem{lagarias} R.~Graham, J.~Lagarias, C.~Mallows, A.~Wilks and C.~Yan, \emph{Apollonian circle packings: number theory,} J. Number Theory \textbf{100} (2003), no. 1, 1--45. 

\bibitem{green-finite-field} B.~J.~Green, \emph{Finite field models in additive combinatorics,} Surveys in Combinatorics 2005, London Math. Soc. Lecture Notes \textbf{327}, 1--27.

\bibitem{green-review} \bysame, \emph{Additive combinatorics: review of the book of Tao and Vu,} Bull. Amer. Math. Soc. (N.S.) \textbf{46} (2009), no. 3, 489--497. 

\bibitem{pfr-blog} \bysame, \emph{The Polynomial Freiman-Ruzsa conjecture,} guest blog at Terry Tao's weblog.

\bibitem{green-partiii} \bysame, notes from a 2009 Cambridge Part III course on Additive Combinatorics, Chapter 2, available on the author's webpage.

\bibitem{green-barbados} \bysame, \emph{Approximate structure in additive combinatorics: Barbados-Radcliffe Lectures,} book in preparation.

\bibitem{green-ruzsa} B.~J.~Green and I.~Z.~Ruzsa, \emph{Freiman's theorem in an arbitrary abelian group,} J. Lond. Math. Soc. (2) \textbf{75} (2007), no. 1, 163--175.

\bibitem{green-tao-equivalence} B.~J.~Green and T.~C.~Tao, \emph{An equivalence between inverse sumset theorems and inverse conjectures for the $U^3$-norm,} preprint available at \texttt{http://arxiv.org/abs/0906.3100}.

\bibitem{gromov} M.~Gromov, \emph{Groups of polynomial growth and expanding maps,} Inst. Hautes \'Etudes Sci. Publ. Math. No. \textbf{53} (1981), 53--73. 

\bibitem{helfgott-1} H.~A.~Helfgott, \emph{Growth and generation in $\SL_2(\Z/p\Z)$,} Ann. of Math. (2) \textbf{167} (2008), no. 2, 601--623. 

\bibitem{helfgott-2} \bysame, \emph{Growth in $\SL_3(\Z/p\Z)$,} to appear in J. European Math. Soc.

\bibitem{hoory-linial-widgerson} S.~Hoory, N.~Linial and A.~Wigderson, \emph{Expander graphs and their applications,} Bull. Amer. Math. Soc. \textbf{43}, no. 4 (2006), 439--561.

\bibitem{hrushovski} E.~Hrushovski, \emph{Stable group theory and approximate subgroups,} preprint available at \texttt{http://arxiv.org/abs/0909.2190}.

\bibitem{iwaniec} H.~Iwaniec, \emph{Unpublished notes on sieve theory.}

\bibitem{iwaniec-kowalski} H.~Iwaniec and E.~Kowalski, \emph{Analytic number theory,} American Mathematical Society Colloquium Publications, \textbf{53}. American Mathematical Society, Providence, RI, 2004. xii+615 pp.

\bibitem{kont-oh} A.~Kontorovich and H.~Oh, \emph{Apollonian circle packings and closed horospheres on hyperbolic 3-manifolds,} preprint available at \texttt{http://arxiv.org/abs/0811.2236}.

\bibitem{lubotzky-survey} A.~Lubotzky, \emph{Cayley graphs: eigenvalues, expanders and random walks,} Surveys in combinatorics, 1995 (Stirling), 155--189, London Math. Soc. Lecture Note Ser., \textbf{218}, Cambridge Univ. Press, Cambridge, 1995. 

\bibitem{lps} A.~Lubotzky, R.~Phillips and P.~Sarnak, \emph{Ramanujan graphs,} Combinatorica \textbf{8} (1988), no. 3, 261--277. 

\bibitem{margulis-expander} G.~A.~Margulis, \emph{Explicit constructions of expanders,} Problemy Pereda\v{c}i Informacii \textbf{9} (1973), no. 4, 71--80.

\bibitem{margulis-combinatorica} \bysame, \emph{Explicit constructions of graphs without short cycles and low density codes,} Combinatorica \textbf{2} (1982), no. 1, 71--78. 

\bibitem{mvw} C.~Matthews, L.~Vaserstein and B.~Weisfeiler, \emph{Congruence properties of Zariski-dense subgroups, I,} Proc. London Math. Soc. \textbf{45} (1984), 514--532.

\bibitem{nori} M.~V.~Nori, \emph{On subgroups of $\GL_n(\F_p)$,} Invent. Math. \textbf{88} (1987), no. 2, 257--275.

\bibitem{razborov} A.~A.~Razborov, \emph{A product theorem in free groups,} preprint available on the author's webpage.

\bibitem{ruzsa-sumset-1} I.~Z.~Ruzsa, \emph{On the cardinality of $A+A$ and $A-A$,} Combinatorics (Proc. Fifth Hungarian Colloq., Keszthely, 1976), Vol. II, pp. 933--938, 
Colloq. Math. Soc. J\'anos Bolyai, 18, North-Holland, Amsterdam-New York, 1978. 

\bibitem{ruz-frei} \bysame, \emph{Generalized arithmetical progressions and sumsets,} Acta Math. Hungar. \textbf{65} (1994), no. 4, 379--388. 

\bibitem{ruzsa-finite-field} \bysame, \emph{An analogue of Freiman's theorem in groups,} Ast\'erisque \textbf{258} (1999), 323--326.

\bibitem{sanders} T.~Sanders, \emph{From polynomial growth to metric balls in polynomial groups,} preprint.

\bibitem{sarnak-expander} P.~Sarnak, \emph{What is$\dots$an expander?} Notices Amer. Math. Soc. \textbf{51} (2004), no. 7, 762--763.

\bibitem{sarnak-survey} \bysame, \emph{Equidistribution and primes,} available on the author's website.

\bibitem{sarnak-letter} \bysame, \emph{Letter to J.~Lagarias,} available on the author's website.

\bibitem{sarnak-xue} P.~Sarnak and X.~X.~Xue, \emph{Bounds for multiplicities of automorphic representations,} Duke Math. J. \textbf{64} (1991), no. 1, 207--227. 

\bibitem{shalom-tao} Y.~Shalom and T.~Tao, \emph{A finitary version of Gromov's polynomial growth theorem,} preprint available at \texttt{http://arxiv.org/abs/0910.4148.}

\bibitem{solymosi-c} J.~Solymosi, \emph{On sum-sets and product-sets of complex numbers,} J. Th\'eor. Nombres Bordeaux \textbf{17} (2005), no. 3, 921--924. 

\bibitem{solymosi-43} \bysame, \emph{An upper bound on the multiplicative energy,} preprint available at \texttt{http://arxiv.org/abs/0806.1040}.

\bibitem{suzuki} M.~Suzuki, \emph{Group theory I,} Springer-Verlag, New York, 1982.

\bibitem{tao-noncommutative} T.~C.~Tao, \emph{Product set estimates for non-commutative groups,} Combinatorica \textbf{28} (2008), no. 5, 547--594. 

\bibitem{tao-general-rings} \bysame, \emph{The sum-product phenomenon in arbitrary rings,} to appear in Contrib. Discrete Math.

\bibitem{tao-solvable} \bysame, \emph{Freiman's theorem for solvable groups,} preprint available at \texttt{http://arxiv.org/abs/}\\
\texttt{0906.3535}.

\bibitem{tao-blog} T.~C.~Tao and others, discussion at \texttt{http://terrytao.wordpress.com/2009/06/21/freimans-}\\
\texttt{theorem-for-solvable-groups/}.

\bibitem{tao-vu-book} T.~C.~Tao and V.~H.~Vu, \emph{Additive Combinatorics,} Cambridge Studies in Advanced Mathematics \textbf{105}, Cambridge University Press 2006.



\end{thebibliography}

\providecommand{\bysame}{\leavevmode\hbox to3em{\hrulefill}\thinspace}

\end{document}